\documentclass[11pt,reqno]{amsart} 
 \usepackage{amsmath,amssymb,amsthm,upref}
 \usepackage{graphicx}
 \usepackage{psfrag}
 \usepackage[numbers,sort&compress,longnamesfirst]{natbib}

\theoremstyle{plain}
\newtheorem{theorem}{Theorem}[section]  
  
\newtheorem{lemma}[theorem]{Lemma}  

\newtheorem*{theorem*}{Theorem}

\theoremstyle{plain}  
\newtheorem{remark}[theorem]{Remark}

\newtheoremstyle{citing}
  {3pt}
  {3pt}
  {\itshape}
  {}
  {\bfseries}
  {.}
  {.5em}
  {\thmnote{#3}}

\theoremstyle{citing}

\numberwithin{equation}{section}
\allowdisplaybreaks[1]

\newlength{\intwidth}
\makeatletter
\DeclareRobustCommand{\cpvint}[2]
    {\mathop{%
       \text{%
         \settowidth{\intwidth}{%
           \ifx\ilimits@\displaylimits
             $\int_{#1}^{#2}$%
           \else
             $\int$%
           \fi}%
         \makebox[0pt][l]{\makebox[\intwidth]{$\text{C}$}}%
         $\int_{#1}^{#2}$}}}
\makeatother

\makeatletter
\DeclareRobustCommand{\cpvintsmall}[2]
    {\mathop{%
       \text{%
         \settowidth{\intwidth}{%
           \ifx\ilimits@\displaylimits
             $\int_{#1}^{#2}$%
           \else
             $\int$%
           \fi}%
         \makebox[0pt][l]{\makebox[\intwidth]{$\text{{\tiny C}}$}}%
         $\int_{#1}^{#2}$}}}
\makeatother

\renewcommand{\a }{\alpha }

\newcommand{\Did }{{\mathcal{D}^{1,2}(\R^3)} } 
\newcommand{\stereo }{\mathcal{S} }

\newcommand{\R}{\mathbb{R}}

\newcommand{\embed}{\hookrightarrow}

\newcommand{\dist}{\text{dist}}

\newcommand{\xpfeil}{\xrightarrow}  
  
\newcommand{\rand}{\partial} 
\newcommand{\where}{\,:\:}  
\newcommand{\iso}{\cong}  
  
\newcommand{\sgn}{\text{sgn}}
  
\newcommand{\supm}{\mathop{\sup}\limits}

\newcommand{\intg}{\mathop{\int}\limits}  
  
\newcommand{\supp}{\mbox{supp}}

\newcommand{\laplace}{\Delta}  
  
\newcommand{\ind}{\text{ind}}  

\newcommand{\di}{\;d}  
  
\newcommand{\nz}{{\mathbb N}}  
\newcommand{\pz}{{\mathbb P}}  
  
\newcommand{\rz}{{\mathbb R}}

\newcommand{\eps}{\varepsilon}  
\renewcommand{\phi}{\varphi} 
\newcommand{\eval}{\vert}

\begin{document}
 
\title[] {The Scalar Curvature Equation on $S^3$}
\author{Matthias Schneider}
\address{Ruprecht-Karls-Universit\"at\\
         Mathematisches Institut\\ 
         Im Neuenheimer Feld 288\\
         69120 Heidelberg, Germany}
\email{mschneid@mathi.uni-heidelberg.de} 
\date{8 April, 2008}  

\begin{abstract}
We give existence results for solutions of the prescribed scalar curvature
equation on $S^3$, when the curvature function is a positive Morse function
and satisfies an index-count condition. 
\end{abstract}

\maketitle

\section{Introduction}
\label{sec:introduction}
Let $S^3$ be the standard sphere with round metric $g_0$ induced by $S^3 = \rand B_1(0) \subset \rz^{4}$.
We study the problem: 
Which functions $K$ on $S^3$ occur as scalar curvature of metrics 
$g$ conformally equivalent to $g_0$? 
Writing $g= \varphi^{4} g_0$ and 
$
k(\theta):= \frac{1}{6}(K(\theta)-6)  
$
this is equivalent to solving for $t=1$ (see \cite{Aubin98})
\begin{align}
\label{eq:eq1}
 -8\laplace_{S^3}\varphi + 6\varphi= 6(1+t k(\theta)) \varphi^{5},\;  \varphi>0 \ \text{ in } S^3.
\end{align}
An obvious necessary condition for the existence of solutions to \eqref{eq:eq1}
is that the function $K$ has to be positive somewhere. Moreover, there 
are the {\em Kazdan-Warner} obstructions \cite{KazdanWarner75,BourguignonEzin87}, which imply 
in particular, that a monotone function of the coordinate function $X_1$
cannot be realized as the scalar curvature of a metric conformal to
$g_0$.\\ 
Numerous studies have been made on equation \eqref{eq:eq1} and its higher dimensional analogue 
and various sufficient conditions
for its solvability have been found (see \cite{AmbAzPer99,YYLi96,YYLi95,ChenLin01,ChenLi01,Bianchi96,AubinBahri97}
and the reference therein), usually under a nondegeneracy assumption on $K$. On $S^3$ a positive function
$K$ is nondegenerate, if
\begin{align}
\label{eq:118}
\tag{nd}
\laplace K(\theta) \neq 0 \text{ if } \nabla K(\theta) =0.   
\end{align}   
For positive Morse functions $K$ on $S^3$ it is shown in 
\cite{SchoenZhang96,BahriCoron91,ChangGurskyYang93} 
that \eqref{eq:eq1} is solvable if $K$ satisfies \eqref{eq:118} and
\begin{align}
\label{eq:34}
d:= -\Big( 1+\sum_{\substack{\nabla K(\theta)=0,\\ \laplace K(\theta)<0}} (-1)^{\ind(K,\theta)}\Big) \neq 0,
\end{align}
where $\ind(K,\theta)$ is the Morse index of $K$ at $\theta$, 
i.e. the number of negative eigenvalues of the Hessian.
For example the simplest possible positive Morse function $K=2+X_1$, where we already know from the Kazdan-Warner
obstructions, that there are no solutions, yields $d=0$, as the only critical point of $K$ with negative Laplacian
is the global maximum with Morse index $3$. Moreover, consider the functions $K_i \in C^{\infty}(S^{3},\rz)$
defined by
\begin{align*}
K_1(X) &:= 2X_1^{2}+6X_2^{2}+7X_3^{2}+8X_4^{2},\\
K_2(X) &:= 3X_1^{2}+6X_2^{2}+7X_3^{2}+8X_4^{2},\\
K_3(X) &:= 4X_1^{2}+6X_2^{2}+7X_3^{2}+8X_4^{2},  
\end{align*}
where $X_i$ for $1\le i\le 4$ is the $i$th coordinate function of $S^3 \subset \rz^4$.
Each $K_i$ is a positive Morse function with critical points given by 
$$\{\pm E_i\in S^{3}\subset \rz^{4} \where  1\le i \le 4\},$$ 
where $\{E_i,1 \le i\le 4\}$ denotes
the standard basis of $\rz^{4}$. The global maximum is attained at $\pm E_4$, the global minimum at $\pm E_1$,
$\pm E_2$ and $\pm E_3$ are saddle points. The sign of the Laplacian, the Morse-index, and $d$ are collected
in Table \ref{tab:degree_k123} below.
\begin{table}[ht]
\begin{tabular}{c||c|c|c}
     & $K_1$ & $K_2$ & $K_3$\\
\hline\hline
$\pm E_1$ & $\laplace_{S^{3}} K_1(\pm E_1)>0$ & $\laplace_{S^{3}} K_2(\pm E_1)>0$ & $\laplace_{S^{3}} K_3(\pm E_1)>0$\\
          & $\ind(K_1,\pm E_1)=0$             & $\ind(K_2,\pm E_1)=0$             & $\ind(K_3,\pm E_1)=0$\\
\hline      
$\pm E_2$ & $\laplace_{S^{3}} K_1(\pm E_2)<0$ & $\laplace_{S^{3}} K_2(\pm E_2)=0$ & $\laplace_{S^{3}} K_3(\pm E_2)>0$\\
          & $\ind(K_1,\pm E_2)=1$             & $\ind(K_2,\pm E_2)=1$             & $\ind(K_3,\pm E_2)=1$\\
\hline      
$\pm E_3$ & $\laplace_{S^{3}} K_1(\pm E_3)<0$ & $\laplace_{S^{3}} K_2(\pm E_3)<0$ & $\laplace_{S^{3}} K_3(\pm E_3)<0$\\
          & $\ind(K_1,\pm E_3)=2$             & $\ind(K_2,\pm E_3)=2$             & $\ind(K_3,\pm E_3)=2$\\
\hline      
$\pm E_4$ & $\laplace_{S^{3}} K_1(\pm E_4)<0$ & $\laplace_{S^{3}} K_2(\pm E_4)<0$ & $\laplace_{S^{3}} K_3(\pm E_4)<0$\\
          & $\ind(K_1,\pm E_4)=3$             & $\ind(K_2,\pm E_4)=3$             & $\ind(K_3,\pm E_4)=3$\\
\hline\hline
d & $1$ & ? &$-1$
\end{tabular} 
\caption{Degree for $K_1$, $K_2$, $K_3$.}
\label{tab:degree_k123}
\end{table}
Thus, (\ref{eq:eq1}) is solvable for $t=1$ and $K \in \{K_1,K_3\}$.
The function $K_2$ does not satisfy the nondegeneracy assumption (\ref{eq:118}) at $E_2$ and the above result 
is not applicable. For the special function $K_2$ a different approach leads to a solution: 
$K_2$ is symmetric with respect to reflections on the sphere $S^{3}$ and the problem may be shifted 
to the projective space $\rz\pz^{3}$. Since $\rz\pz^{3}$ is not conformal to $S^{3}$ the result of
\citet{MR856845} yields a solution on $\rz\pz^{3}$ that may be shifted back to obtain a solution for
$K_2$ on $S^{3}$. But, the argument breaks down for any nonsymmetric perturbation of $K_2$. 
We are interested exactly in this case,
when the nondegeneracy assumption \eqref{eq:118} is not satisfied, and we shall give the required general existence
result.\\     
In the following, unless otherwise stated, 
we will assume that $K=6(1+k) \in C^5(S^3)$ is positive.
To give our main results we need the following notation. 
We denote by $\stereo_{\theta}(\cdot)$ stereographic coordinates 
centered at some point $\theta \in S^3$,
i.e. $\stereo_{\theta}(0)=\theta$.
We write $k_\theta= k \circ \stereo_\theta$ and
for a critical point $\theta$ of $k$ with $D^2k_{\theta}(0)$ invertible we let 
\begin{align*}
a_0(\theta)&:= \cpvint {\rz^3}{} \Big(k_{\theta}(x)-T_{k_\theta,0}^2(x)\Big)|x|^{-6}, \\
a_1(\theta)&:=  \laplace^2k_{\theta}(0)
+\nabla(\laplace k_{\theta}(0)) \cdot \big(D^2k_{\theta}(0)\big)^{-1}\nabla(\laplace k_{\theta}(0)),\\
a_2(\theta)&:= k_{\theta}(0)a_1(\theta)
- \frac{15}{8\pi}\int_{\rand B_1(0)}\big|D^2k_{\theta}(0)(x)^2\big|^2,
\end{align*}  
where all differentiations are done in $\rz^3$, the $m$th Taylor polynomial of $k_\theta$ in $y$
is abbreviated by  
\begin{align*}
T_{k_\theta,y}^m(x) := \sum_{\ell=0}^m\frac{1}{\ell!}D^\ell k_\theta(y)(x-y)^\ell,
\end{align*}
and $\cpvintsmall{}{}$ is the Cauchy principal value of the integral,
\begin{align*}
\cpvint{\rz^3}{} f(x) := \lim_{r \to 0} \int_{\rz^3\setminus B_r(0)}f(x). 
\end{align*}
The value $a_0(\theta)$ is well defined because of the cancellation due to symmetry. 
For instance expanding
$T_{k_\theta,0}^{m}$ in spherical harmonics we get
\begin{align*}
\int_{\rand B_1(0)} T_{k_\theta,0}^{m}(x)\, dS =
\begin{cases}
0 &\text{if } m \text{ is odd,}\\
\frac{2 \pi}{3} \laplace k_\theta(0)&\text{if }m=2.   
\end{cases}
\end{align*}
The value $a_0(\theta)$ will be of interest only in points 
where (\ref{eq:118}) is not satisfied, that is when
$\nabla k_\theta(0)$ and $\laplace k_\theta(0)$ vanish simultaneously. 
In this case $a_0(\theta)$ is given by
\begin{align*}
a_0(\theta) = \cpvint {\rz^3}{} \Big(k_{\theta}(x)-k_\theta(0)\Big)|x|^{-6},  
\end{align*}
and measures, weighted by $|x|^{-6}$, the difference between $k_{\theta}$ and $k_\theta(0)$.\\ 
Denote by $\text{Crit}(k)$, $M$, and $T$ the sets,
\begin{align*}
\text{Crit}(k) &:= \big\{\theta \in S^3 \where \nabla k(\theta)=0 \big\},\\
M &:= \big\{\theta \in \text{Crit}(k)\where \laplace k_\theta(0)=a_0(\theta)= 0, 
\,\text{and } a_2(\theta)\neq 0\big\},\\
T &:= \{-a_1(\theta)/a_2(\theta) \where \theta \in M\},
\end{align*}
\begin{theorem}
\label{t:main}
Suppose $1+k\in C^5(S^3)$ is a positive Morse function. 
Then \eqref{eq:eq1} is solvable for $t \in (0,1]\setminus T$, if 
\begin{align}
\label{eq:11}
0 \neq d(t)= -\Big(1 + \sum_{\theta \in \text{Crit}_-(k,t)} (-1)^{\ind(k,\theta)}\Big),
\end{align}
where
\begin{align*}
\text{Crit}_-(k,t) := \Big\{&\theta \in S^3 \where \nabla k(\theta)=0 
\text{ and }\\ 
&\lim_{\mu \to 0^+} \sgn\Big(\laplace k(\theta)+ a_0(\theta)\mu -\big(a_1(\theta)+ta_2(\theta)\big)\mu^2\Big)= -1\Big\}. 
\end{align*}
The number $d(t)$ is the Leray-Schauder degree of the problem \eqref{eq:eq1}.
\end{theorem}
We note that set of critical points of $K$ and $k$ are equal and for any $\theta \in \text{Crit}(k)$ we have
\begin{align*}
\sgn(\laplace_{S^{3}}K(\theta)) &= \sgn(\laplace_{S^{3}}k(\theta))= \sgn(\laplace_{\rz^{3}}k_\theta(0)),\\
\ind(K,\theta) &= \ind(k,\theta)= \ind(k_\theta,0).  
\end{align*}
Hence, the nondegeneracy condition (\ref{eq:118}) implies that the set $M$ is empty and the formula
in (\ref{eq:11}) gives exactly the index-count condition in (\ref{eq:34}).
In contrast to (\ref{eq:34}) the Leray-Schauder degree now depends on $t$ and may change as $t$ crosses
some value in $T$. Indeed for any 
\begin{align*}
t_*=-\frac{a_1(\theta)}{a_2(\theta)} \in T \cap (0,1]  
\end{align*}
there is a ``blow-up curve'' $(t(s),\phi(s))$ such that 
\begin{align*}
\lim_{s \to 0}t(s)=t_*,\;
\lim_{s \to 0} \|\phi(s)\|_{L^{\infty}(B_\eps(\theta))}=+\infty \text{  for all }\eps>0,
\end{align*}
and $\phi(s)$ solves (\ref{eq:eq1}) with $t=t(s)$ (see \cite{schneider04_2} and Figure \ref{fig:blow_up_curves} below).
\begin{figure}[ht]
\resizebox{10cm}{!}{\includegraphics{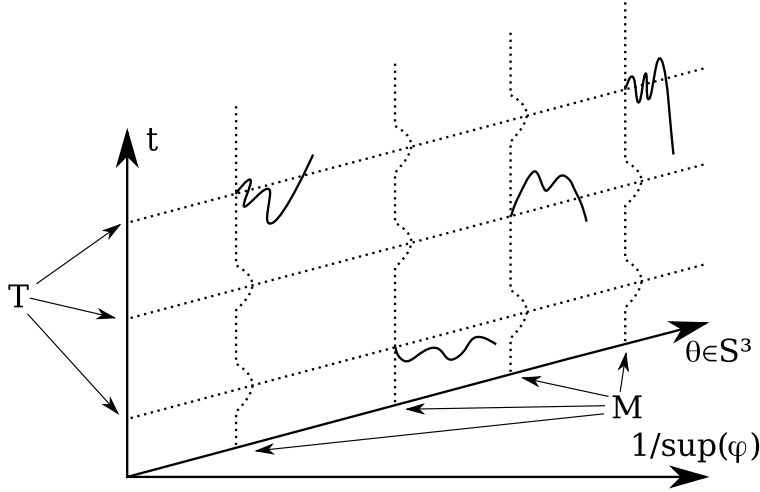}}
\caption{Blow up curves}
\label{fig:blow_up_curves}
\end{figure}\\
An inspection of the proof of Theorem \ref{t:main} shows that the result remains valid, when
$k$ is only in $C^4(S^3)$. We state Theorem \ref{t:main} for functions $k \in C^5(S^3)$, 
because we use the analysis in \cite{schneider04_2,schneider04}, which
is done in this setting.\\
To illustrate our results
we will apply Theorem \ref{t:main} when $K$ equals $K_i$ for some $i \in \{1,2,3\}$. 
For $i\in \{1,3\}$ the set $M$
is empty, as the Laplacian does not vanish at any critical point, 
$d(\cdot)$ is independent of $t\neq 0$ and given by
(\ref{eq:34}). Concerning $K_2$, the critical points with vanishing Laplacian are $\{\pm E_2\}$ 
and we need to compute $a_j(\pm E_2)$ for
$j=0,1,2$ and the function
\begin{align*}
k=k_2:= \frac{1}{6}(K_2-6)= \frac{1}{2} X_1^{2}+X_2^{2}+\frac{7}{6}X_3^{2}+\frac{4}{3}X_4^{2}-1.  
\end{align*}
A straightforward computation (see \cite{schneiderhabil}) shows
\begin{align*}
a_0(\pm E_2)=0,\; a_1(\pm E_2)=0,\; a_2(\pm E_2)= -\frac{224}{9}.  
\end{align*}
Hence, $M=\{\pm E_2\} \subset S^3$, $T=\{0\}$, and 
\begin{align*}
d(t)=
\begin{cases}
-1 &\text{if } t>0,\\
 1 &\text{if } t<0.  
\end{cases}
\end{align*}
Thus, we may replace the question mark in Table \ref{tab:degree_k123} by $-1$.
Moreover, for $0\neq h\in
C_c^{\infty}(S^{3}\setminus\{\pm E_2\},\rz_{\ge 0})$ we consider $k_2 \pm s h$, where $s$ is a small positive
parameter. Since
\begin{align*}
\int_{\rz^3} h_{\pm E_2}(x)|x|^{-6}\, dx >0,  
\end{align*}
the sets $M$ and $T$ are empty for $k=k_2\pm s h$ and $0<s<<1$, the degree for $t\neq 0$ is given by
\begin{align*}
d(t)=-1 \text{ for } k=k_2 + s h, \; d(t)=1 \text{ for } k=k_2 - s h.  
\end{align*}
Furthermore, we consider for $0<s<<1$
\begin{align*}
k = k_2 +s\big(7(1-X_2^{2})^{2}-20(1-X_2^{2})^{3}\big)  
\end{align*}
For small positive $s$ the set of critical points of $K$ is given by $\{\pm E_i\}$ with vanishing Laplacian only at
$\pm E_2$, $a_0(\pm E_2)=0$, and
\begin{align*}
a_1(\pm E_2)= 13440s,\, a_2(\pm E_2)= -\frac{224}{9}.   
\end{align*}
Thus, $M=\{\pm E_2\}$, $T=\{540s\}$, and for $t\neq 0$
\begin{align*}
d(t)= 
\begin{cases}
-1 &\text{if }t>540s,\\
1&\text{if }t<540s. 
\end{cases}
\end{align*}
The change of the degree is due to the two blow-up curves $r \mapsto (t^{\pm}(r),\phi^{\pm}(r))$,
where $t^{\pm}(r) \to 540s$ and $\phi^{\pm}(r)$ concentrates at $\pm E_2$ as $r \to 0$.
It is interesting to note that, although $K$ is even in this case, the solutions on the blow-up curve
are not even as they concentrate in a single point.\\
To prove our main result we embed our problem into a two dimensional family of problems.
We choose $h\in C^\infty(S^3,[0,\infty))$ such that 
\begin{align}
\label{eq:4}
\supp(h) \cap \text{Crit}(k) = \emptyset.  
\end{align}
We fix $0<t_0 \in (0,1]\setminus T$ and consider for $s \ge 0$
\begin{align}
\label{eq:10}
-8\laplace_{S^3}\varphi + 6\varphi= 6\Big(1+t_0 \big(k(\theta)+sh(\theta)\big)\Big) \varphi^{5},\;  
\varphi>0 \ \text{ in } S^3.  
\end{align}
Analogously as above, we define $a_j(\theta,s)$ for $j=0,1,2$ and $M_s$ by replacing $k$ by $k+sh$
in the definition of $a_j(\theta)$ and $M$. We obtain for $\theta \notin \supp(h)$ 
\begin{align*}
a_0(\theta,s)= a_0(\theta) +s \int_{\rz^3} h_\theta(x)|x|^{-6},\;
a_1(\theta,s)=  a_1(\theta),\; a_2(\theta,s)= a_2(\theta).
\end{align*}
From (\ref{eq:4}) there is $s_0>0$ such that for $0 \le s \le s_0$: 
\begin{itemize}
\item $\text{Crit}(k)=\text{Crit}(k+sh)$,
\item $k+sh$ is a Morse function,
\item $a_0(\theta)\cdot a_0(\theta,s)>0$, if $\nabla k(\theta) =0$ and $a_0(\theta) \neq 0$.  
\end{itemize}
The main reason for introducing the perturbation $h$ is that the sets $M_s$ are empty, because
\begin{align*}
a_0(\theta,s) \neq 0, \text{ if } \nabla k(\theta) =0.  
\end{align*}
By standard elliptic regularity the operator $L_{s}$, defined by
\begin{align*}
L_s: \: \phi \mapsto (-8\laplace_{S^3}+6)^{-1}\bigg(6\Big(1+t_0 \big(k(\theta)+sh(\theta)\big)\Big)\phi^5\bigg),  
\end{align*}
is compact from $C^2(S^3)$ into $C^2(S^3)$. From the apriori estimates in \cite{schneider04_2}
, as $t_0 \notin T$, there is $C_{t_0}>0$ such that all positive solution
to (\ref{eq:10}) with $s=0$ lie in $\mathcal{B}_{C_{t_0}}$,
\begin{align*}
\mathcal{B}_{C} := \{\phi \in C^2(S^3) \where \|\phi\|_{C^2(S^3)}<C 
\text{ and }C^{-1}< \phi\}. 
\end{align*}
Moreover, as $\text{Crit}(k+sh)$ does not change when $s$ moves from $0$ to $s_0$, we may apply
Theorem 7.1 in \cite{schneider04_2}. Thus, for any $0<\delta<s_0$ there is $C_\delta>0$ such that 
all positive solution
to (\ref{eq:10}) with $s \in [\delta,s_0]$ lie in $\mathcal{B}_{C_{\delta}}$. 
The Leray-Schauder degree 
$\deg(Id-L_s,\mathcal{B}_{C_{\delta}},0)$, which is well-defined and independent of $s \in [\delta,s_0]$ by the apriori 
estimates, is computed in \cite{schneider04} and equals
\begin{align}
\label{eq:16}
\deg(Id-L_s,\mathcal{B}_{C_{\delta}},0)= -\Big(1+ \sum_{\theta \in Crit_-(k+sh)} (-1)^{\ind(k,\theta)}\Big), 
\end{align}
where the set $Crit_-(k+sh)$ is given by
\begin{align*}
Crit_-(k+sh) &:= \Big\{ \theta \in \text{Crit}(k) \where 
\lim_{\mu \to 0^+} \sgn\Big(\laplace k(\theta)+ a_0(\theta,s)\mu\Big)= -1\Big\}.
\end{align*}
As $h\ge 0$, we have for $\theta \in \text{Crit}(k)$ 
that
$a_0(\theta,s)<0$ if and only if $a_0(\theta)<0$.  
hence
\begin{align*}
Crit_-(k+sh)
&= \big\{\theta \in \text{Crit}(k) \where \laplace k(\theta)<0 \text{ or } 
\big(\laplace k(\theta)=0 \text{ and } a_0(\theta)<0\big)\big\}.
\end{align*}
The constant
$C_{\delta}$ in \cite{schneider04_2} depends on 
\begin{align*}
\sup_{s \in [\delta,s_0]} 
\{|a_0(\theta,s)|^{-1} \where \nabla k(\theta)=0,\,\laplace k_{\theta}(0) =0 \text{ and } a_0(\theta,s)\neq 0\}.
\end{align*}
Consequently, we cannot assume that $C_\delta$ remains bounded as $\delta \to 0$. Indeed, we shall show
that as $s$ moves to $0$ the family of solutions splits into solutions, that remain uniformly bounded 
as $s \to 0^+$ and converge to solutions of (\ref{eq:10}) with $s=0$, and solutions that blow up as $s \to 0^+$. 
When $s$ moves to $0^+$ the total degree,
which is computed in (\ref{eq:16}), is given by the sum of two degree's, the degree of the ``bounded solutions'',
that we are interested in, and the degree of the ``blow-up solutions''.
We will compute the degree of the solutions, that blow up when $s \to 0^+$, as a sum of local degree's.
Subtracting the result from (\ref{eq:16}) leads to the formula in (\ref{eq:11}).  
 
\section{Preliminaries}
\label{s:preliminiaries}
For fixed $\theta \in S^{3}$ in stereographic coordinates $\stereo_{\theta}(\cdot)$ 
equation (\ref{eq:10})
is equivalent to
\begin{align}
\label{eq:eq2}
-\laplace u = \big(1+t_0 (k_\theta(x)+s h_\theta(x))\big) u^{5} \text{ in }\rz^3, \; u>0. 
\end{align}
where $h_\theta = h \circ\stereo_{\theta}$ and 
\begin{align}
\label{eq:116}
u(x) = {\mathcal R}_\theta(\phi)(x):= 3^{\frac{1}{4}}(1+|x|^2)^{-\frac{1}{2}} \phi \circ \stereo_{\theta}(x).  
\end{align}
The transformation (\ref{eq:116}) gives rise to a Hilbert space isomorphism between
$H^{1,2}(S^3)$ and $\Did$, the closure of $C_c^\infty(\rz^3)$ with respect to 
\begin{align*}
\|u\|^2:= \int_{\rz^3} |\nabla u|^2 =\langle u,u\rangle.   
\end{align*}
Due to elliptic regularity (see \cite{BrezisKato79,Luckhaus79}) and
Harnack's inequality it is enough to find a weak nonnegative solution of (\ref{eq:10}) in $H^{1,2}(S^3)$, or of
the equivalent equation in $\Did$. Although we take advantage of both formulations, 
we mainly  consider (\ref{eq:eq2}) to analyze the blow-up behavior and to compute local degrees.
Weak solutions to (\ref{eq:eq2}) correspond to critical points of $f_{t_0,s}: \Did \to \rz$
\begin{align*}
f_{t_0,s}(u) :=  \int_{\rz^{3}} \frac12 |\nabla u|^{2} 
-\frac{1}{6}\Big(1+t_0 \big(k_\theta(x)+s h_\theta(x)\big)\Big) u^{6}\, dx.
\end{align*}
We denote by $f_0$ the unperturbed functional with $t_0=s=0$. 
The positive solutions of \eqref{eq:eq2} for $t_0=s=0$, i.e. the positive critical points of $f_0$, 
are completely known 
(see \cite{CaffarelliGidasSpruck89,GidasNiNirenberg81,ChenLi91})
and given by a noncompact manifold  
\begin{align*}
Z:= \big\{z_{\mu,y}(x):= \mu^{-\frac{1}{2}}3^{\frac{1}{4}} 
\Big(1+\big|\frac{x-y}{\mu}\big|^2\Big)^{-\frac{1}{2}} \where y \in \rz^3, \mu>0\big\},
\end{align*}
We state some properties of the critical manifold $Z$ and $f_0$ (see \cite{schneider04_2,AmbAzPer99} for details).
We define for $\mu>0$ and $y \in \rz^3$ the 
maps ${\mathcal U}_{\mu},\,{\mathcal T}_y: \Did \to \Did$ by
\begin{align*}
{\mathcal U}_{\mu}(u):= \mu^{-\frac{1}{2}} u\Big(\frac{\cdot}{\mu}\Big) \text{ and }
{\mathcal T}_{y}(u):= u(\cdot-y).  
\end{align*}
With this notation the critical manifold $Z$ is given by
\begin{align*}
Z = \{ z_{\mu,y}= {\mathcal T}_y \circ {\mathcal U}_\mu (z_{1,0}) \where y \in \rz^3,\, \mu>0\}.
\end{align*}
The dilation ${\mathcal U}_{\mu}$ and the translation ${\mathcal T}_{y}$ are automorphisms 
of $\Did$ and for every $\mu>0$, $y \in \rz^3$, and  $v \in \Did$ 
\begin{align}
\label{eq:7}
\begin{split}
({\mathcal U}_{\mu})^{-1} = ({\mathcal U}_{\mu})^{t} = U_{\mu^{-1}},\
({\mathcal T}_{y})^{-1} = ({\mathcal T}_{y})^{t} = {\mathcal T}_{-y},\\ 
f_0 = f_0 \circ {\mathcal U}_{\mu} = f_0 \circ {\mathcal T}_{y}, \text{ and }\\
f_0''(v) = 
({\mathcal T}_{y}\circ {\mathcal U}_{\mu})^{-1} \circ 
f_0''({\mathcal T}_{y}\circ {\mathcal U}_{\mu}(v)) \circ ({\mathcal T}_{y}\circ {\mathcal U}_{\mu}),
\end{split}  
\end{align}
where $(\cdot)^{t}$ denotes the adjoint. 
The tangent space $T_{z_{\mu,y}}Z$ at a point $z_{\mu,y} \in Z$ is spanned by $4$ orthonormal functions,
\begin{align*}
T_{z_{\mu,y}}Z = \langle (\dot{\xi}_{\mu,y})_i \where i=0 \dots 3 \rangle,\\
(\dot{\xi}_{\mu,y})_i :=
\begin{cases}
\|\frac{d}{d\mu} z_{\mu,y}\|^{-1} \frac{d}{d\mu} z_{\mu,y} &\text{, if }i=0,\\
\|\frac{d}{d y_i} z_{\mu,y}\|^{-1} \frac{d}{d y_i} z_{\mu,y} &\text{, if }1\le i\le 3.  
\end{cases}
\end{align*}
The maps ${\mathcal U}_{\mu}$ and ${\mathcal T}_{y}$ are isomorphism of the tangent spaces, and moreover 
\begin{align}
\label{eq:2}
\begin{split}
({\dot\xi}_{\mu,y})_i= {\mathcal T}_y \circ {\mathcal U}_\mu \big(({\dot\xi}_{1,0})_i\big),\\
{\mathcal T}_{y}\circ{\mathcal U}_{\mu}:\: 
(T_zZ)^{\perp} \xpfeil{\iso} 
(T_{{\mathcal T}_{y}\circ{\mathcal U}_{\mu}(z)}Z)^{\perp}.   
\end{split}
\end{align}
We consider $f_{t_0,s}'(u)$ as an element of $\Did$ and $f_{t_0,s}''(u)$ as a map in
$\mathcal{L}(\Did)$. With this identification $f_{t_0,s}''(u)$ 
is a self-adjoint, 
compact perturbation of the identity map in $\Did$. The spectrum $\sigma(f_0''(z_{\mu,y}))$
consists of point-spectrum accumulating at $1$
and is computed together with the eigenspaces in \cite{schneider04_2}. 
Since $Z$ is a manifold of critical points of $f_0'$, the tangent space $T_z Z$ at a point $z \in Z$ is contained 
in the kernel $N(f_0''(z))$ of $f_0''(z)$, knowing the eigenspaces we see 
\begin{align}
\label{eq:15}
T_z Z=  N(f_0''(z)) \quad \text{ for all } z \in Z. 
\end{align}
If \eqref{eq:15} holds the critical manifold $Z$ is called nondegenerate (see \cite{AmBa2}).
The operator $f_0''(z)$ maps the space $\Did$ into
$T_{z}Z^\perp$ and is invertible in ${\mathcal{L}}(T_{z}Z^\perp)$. 
From \eqref{eq:7} and \eqref{eq:2}, we obtain in this case 
\begin{align}
\label{eq:3}
\|(f_0''(z_{1,0}))^{-1}\|_{{\mathcal{L}}(T_{z_{1,0}}Z^\perp)} = 
\|(f_0''(z))^{-1}\|_{{\mathcal{L}}(T_{z}Z^\perp)} \quad \forall z \in Z.  
\end{align}
Moreover, $T_{z_{\mu,y}}Z^\perp$ splits orthogonally into (see \cite{schneider04_2})
\begin{align}
\label{eq:5}
T_{z_{\mu,y}}Z^\perp = \langle z_{\mu,y}\rangle \oplus^{\perp }
\langle \Phi^{\mu,y}_{i,j,l}:\, i,j\in \nz_0,\, 2\le i+j\le n ,\, 1\le l\le  c_i \rangle,  
\end{align}
where $\Phi^{\mu,y}_{i,j,l}$ are eigenfunctions of $f_0''(z_{\mu,y})$ with positive eigenvalue
\begin{align*}
\lambda_{i,j} = 1-\frac{15}{(4+2(i+j-1))^{2}-1}.  
\end{align*}
The dimension of the eigenspace corresponding to $\lambda_{i.j}$ is denoted by $c_i$.
The functions $\Phi^{\mu,y}_{i,j,l}$ are smooth and given in terms of Jacobi polynomials and spherical harmonics. 
The operator $f_0''(z_{\mu,y})$ has precisely one negative 
eigenvalue $-4$ with one-dimensional eigenspace $\langle z_{\mu,y}\rangle$. 

\section{The blow up analysis}
\label{s:blow_up}
Based on the results in \cite{SchoenZhang96,YYLi95} we have the following lemma (see \cite[Cor 3.2]{schneider04_2})
\begin{lemma}
\label{l:blow_up}  
Suppose $1+t_0k \in C^1(S^3)$ is positive and $h \in C^\infty(S^3)$ is nonnegative.  
If $(s_i,\varphi_i) \in [0,1] \times C^2(S^3)$ solve (\ref{eq:10}) with $s=s_i$,
then after passing to a subsequence  either $(\varphi_i)$
is uniformly bounded in $L^\infty(S^3)$ (and hence in $C^{2,\alpha}(S^3)$ by standard elliptic regularity) 
or there exist $\theta \in S^3$ and sequences $(\mu_i)\in (0,\infty)$, $(y_i) \in \rz^3$ satisfying
$\lim_{i \to \infty} \mu_i= 0$ and $\lim_{i \to \infty} y_i= 0$,       
such that (in stereographic coordinates $\stereo_\theta(\cdot)$)
\begin{align*}
{\mathcal R}_\theta(\phi_i) - \Big(1+t_0\big(k_\theta(y_i)+s_ih_\theta(y_i)\big)\Big)^{-\frac14} z_{\mu_i,y_i} 
\text{ is orthogonal to } T_{z_{\mu_i,y_i}}Z,\\
\|{\mathcal R}_\theta(\phi_i) 
- \Big(1+t_0\big(k_\theta(y_i)+s_ih_\theta(y_i)\big)\Big)^{-\frac14} z_{\mu_i,y_i}\|_{\Did} = o(1).
\end{align*}
\end{lemma}

\section{The finite dimensional reduction}
\label{s:expansion}
For the rest of the paper, unless otherwise indicated, integration
extends over $\rz^3$ and is done with respect to the variable $x$.
\begin{lemma}
\label{p:implicit} 
Suppose $1+k \in C^5(S^3)$ is a positive Morse function, $t_0 \in (0,1]\setminus T$,
$h \in C^\infty(S^3)$ satisfies (\ref{eq:4}), and $\theta \in S^3$. 
Then there exist $s_0= s_0(t_0,k,h)>0$, 
$\mu_0= \mu_0(t_0,k,h)>0$ and two functions $w:\Omega \to \Did$ and $\vec\alpha:\Omega \to \rz^{4}$ depending on
$t_0$, $k$, $h$, and $\theta$, where
\begin{align*}
\Omega:= \{(s,\mu,y) \in (-s_0,+s_0) \times (0,\mu_0)\times \rz^3\}
\end{align*}
such that for any $(s,\mu,y) \in \Omega$
\begin{align}
\label{eq:12}
w(s,\mu,y)\ \text{ is orthogonal to }\ T_{z_{\mu,y}}Z\\
\label{eq:13}
f_{t_0,s}'\big(z_{\mu,y}+w(s,\mu,y)\big)= \vec \alpha(s,\mu,y) \cdot \dot{\xi}_{\mu,y} \in T_{z_{\mu,y}}Z\\
\label{eq:9}
\|w(s,\mu,y) -  w_0(s,\mu,y)\|+ \|\vec{\alpha}(s,\mu,y)\| < \rho_0,
\end{align}
where $\{({\dot\xi}_{\mu,y})_i \where i=0 \dots 3\}$ denotes the basis of $T_{z_{\mu,y}}Z$ 
given in \eqref{eq:11} and
\begin{align*}
w_0(s,\mu,y) := \Big(\big(1+t_0(k_\theta(y)+s h_\theta(y))\big)^{-\frac{1}{4}}-1\Big)z_{\mu,y}.
\end{align*}
The functions $w$ and $\vec \alpha$ are of class $C^{2}$ 
and unique in the sense that if $(v,\vec \beta)$ satisfies
\eqref{eq:12}-\eqref{eq:9} for some $(s,\mu,y)\in \Omega$ then $(v,\vec \beta)$ is given by
$(w(s,\mu,y),\vec \alpha(s,\mu,y))$.\\
Moreover, we have we have as $\mu \to 0$
\begin{align*}
|\vec{\alpha}(&s,\mu,y)-\sum_{j=1}^4\vec{\alpha}_{j}(s,\mu,y)| 
= O(\mu^{4+\frac14}+ \mu^2|\nabla k_\theta(y)+s\nabla h_\theta(y)|^2)\\ 
&+O\big(\mu^3|\nabla k_\theta(y)+s\nabla h_\theta(y)|
+ \mu^4 |\laplace k_\theta(y)+s\laplace h_\theta(y)|\big),    
\end{align*}
where $\a_1$, $\a_2$ are given by
\begin{align*}
\vec{\alpha}_1(s,\mu,y)&:= -\mu \Big(1+t_0 \big(k_\theta(y)+sh_\theta(y)\big)\Big)^{-\frac{5}{4}}\,
\frac{t_0\pi}{3^{\frac14}\sqrt{5}}
\binom{0}{\nabla k_\theta(y)+s\nabla h_\theta(y)},\\
\vec{\alpha}_{2}(s,\mu,y) &:= - \mu^2 \Big(1+t_0 \big(k_\theta(y)+sh_\theta(y)\big)\Big)^{-\frac{5}{4}}
\frac{t_0\pi}{3^{\frac14}\sqrt{5}} \binom{\laplace \big(k_\theta+s h_\theta\big)(y)}{\vec 0},
\end{align*}
for $1\le i \le 3$
\begin{align*}
\vec{\alpha}_{3}(s,\mu,y)_i &:= - \mu^3 \Big(1+t_0 \big(k_\theta(y)+sh_\theta(y)\big)\Big)^{-\frac{5}{4}}\,
\frac{t_0\pi}{3^{\frac14}2\sqrt{5}}
 \frac{\rand}{\rand x_i}\laplace \big(k_\theta +s h_\theta\big)(y),\\
\vec{\alpha}_{4}(s,\mu,y)_i &:= - \mu^4  \Big(1+t_0\big(k_\theta+sh_\theta\big)(y)\Big)^{-\frac{5}{4}}
\frac{t_0 3^\frac34 8}{\pi \sqrt{5}}\\
&\qquad   
\cpvint{}{} \Big((k_\theta+s h_\theta)(x+y)-  T^3_{(k_\theta+sh_\theta)(\cdot +y),0}(x)\Big)\frac{x_i}{|x|^{8}},
\end{align*}
and
\begin{align*}
\vec{\alpha}_{3}(s,\mu,y)_0 &:= - \mu^3 \Big(1+t_0\big(k_\theta+sh_\theta\big)(y)\Big)^{-\frac{5}{4}} 
\frac{t_0 3^\frac34 4}{\pi \sqrt{5}}\\
&\qquad 
\cpvint{}{} \big((k_\theta+sh_\theta)(x+y)- T_{(k_\theta+sh_\theta)(\cdot+y),0}^2(x)\big)\frac{1}{|x|^{6}},\\
\vec{\alpha}_{4}(s,\mu,y)_0 &:=  \mu^4  \Big(1+t_0\big(k_\theta+sh_\theta\big)(y)\Big)^{-\frac{5}{4}} 
\frac{t_0 3^\frac34 \pi \sqrt{5}}{30} \laplace^2 \big(k_\theta+sh_\theta\big)(y)+\\ 
&\frac{-t_0^2 \mu^4\, 3^\frac34 \sqrt{5}}{16\Big(1+t_0\big(k_\theta+sh_\theta\big)(y)\Big)^{\frac94}}
\intg_{\rand B_1(0)} \big|D^2\big(k_\theta+sh_\theta\big)(y)(x)^2\big|^2 \di Sx.
\end{align*}
\end{lemma}
Replacing $k$ by $k+sh$
the existence part, uniqueness, and the asymptotic estimates as $\mu \to 0$ follow directly from Lemmas 4.2-4.7 in
\cite{schneider04_2}. It only remains to show the $C^2$-dependence on $s$, which we omit, since it is analogous to the
proof given in \cite{schneider04_2}.\\
Concerning the derivatives of $\vec{\a}$ with respect to $\mu$ and $y$ we may apply the results in
\cite[Lem. 5.1]{schneider04_2} and \cite[Lem A.4-A.5]{schneider04} to obtain the following two lemmas.
\begin{lemma}
\label{l:derivative_alpha_y}
Under the assumptions of Lemma \ref{p:implicit} we have for all
$(s,\mu,y) \in \Omega$ and $1 \le i,j \le 3$
\begin{align*}
\frac{\rand \alpha(s,\mu,y)_i}{\rand y_j} 
&= - \frac{t_0 \mu \pi}{3^\frac14 \sqrt{5}}\Big(1+t_0 \big(k_\theta+sh_\theta\big)(y)\Big)^{-\frac54} 
\frac{\rand^2 (k_\theta+sh_\theta)(y)}{\rand x_i \rand x_j}\\ 
&\quad +O\big(|\nabla \big(k_\theta+sh_\theta\big)(y)|^2\mu+\mu^{2+\frac14}\big).
\end{align*}
\begin{align*}
\frac{\rand \alpha(s,\mu,y)_0}{\rand y_j} 
&=-\frac{t_0 \mu^2 \pi}{3^\frac14 \sqrt{5}} \Big(1+t_0\big(k_\theta+sh_\theta\big)(y)\Big)^{-\frac54}
\frac{\rand}{\rand x_j} \laplace \big(k_\theta+sh_\theta\big)(y)\\ 
&\quad + O\big(|\nabla (k_\theta+sh_\theta)(y)|^2\mu+ \mu^{2+\frac14}\big),  
\end{align*}
\end{lemma}
 
\begin{lemma}
\label{l:derivative_alpha_mu}
Under the assumptions of Lemma \ref{p:implicit} we have for all
$(s,\mu,y) \in \Omega$ and $1 \le i \le 3$
\begin{align*}
\frac{\rand \alpha(s,\mu,y)_i}{\rand \mu} 
&= \sum_{j=1}^3 \frac{\rand \alpha_j(s,\mu,y)_i}{\rand \mu} 
+O\big((|\nabla k_\theta(y)|^2+|\nabla h_\theta(y)|^2)\mu+\mu^{3}\big),
\end{align*}
\begin{align*}
\frac{\rand \alpha(s,\mu,y)_0}{\rand \mu} 
&= \sum_{j=2}^4 \frac{\rand \alpha_j(s,\mu,y)_0}{\rand \mu} 
 +O\big((|\laplace k_\theta(y)|+|\laplace h_\theta(y)|)\mu^3+\mu^{3+\frac14}\big)\\
&\quad+O\big((|\nabla k_\theta(y)|^2+|\nabla h_\theta(y)|^2)\mu+(|\nabla k_\theta(y)|+|\nabla h_\theta(y)|)\mu^2\big).
\end{align*}
\end{lemma}
In order to compute the derivative of $\vec{\a}$ with respect to $s$ one has to mimic the lengthy
calculation of the $t$-derivative in \cite[Lem 5.2-5.3]{schneider04_2}. 
We will again just state the 
result and refer to \cite{schneiderhabil} for details. 
This will be the last point where we are less precise concerning the $s$-dependence.   
\begin{lemma}
\label{l:derivative_alpha_s}
Under the assumptions of Lemma \ref{p:implicit} we have for all
$(s,\mu,y) \in \Omega$ and $1 \le i \le 3$
\begin{align*}
\frac{\rand \alpha(s,\mu,y)_i}{\rand s} 
&= \sum_{j=1}^3 \frac{\rand \alpha_j(s,\mu,y)_i}{\rand s} 
+O\big((|\nabla k_\theta(y)|^2+|\nabla h_\theta(y)|^2)\mu^2+\mu^{4}\big),
\end{align*}
\begin{align*}
\frac{\rand \alpha(s,\mu,y)_0}{\rand s} 
&= \sum_{j=2}^4 \frac{\rand \alpha_j(s,\mu,y)_0}{\rand s} 
 +O\big((|\laplace k_\theta(y)|+|\laplace h_\theta(y)|)\mu^4+\mu^{4+\frac14}\big)\\
&\quad+O\big((|\nabla k_\theta(y)|^2+|\nabla h_\theta(y)|^2)\mu^2+(|\nabla k_\theta(y)|+|\nabla h_\theta(y)|)\mu^3\Big).
\end{align*}
\end{lemma}

\begin{lemma}
\label{l:implicit:curve}
Under the assumptions of Lemma \ref{p:implicit} suppose 
\begin{align*}
\nabla k_\theta(0)=0 \text{ and }  \laplace k_\theta(0)=0.
\end{align*}
Consider the function $\hat{\alpha}:\Omega \to \rz^3$, defined by
\begin{align*}
\hat{\alpha}(s,\mu,y) := \frac{3^\frac14 \sqrt{5}}{t_0 \mu \pi} 
\big(1+t_0 k(\theta)\big)^{\frac54}(\vec{\alpha}(s,\mu,y)_{1},\dots,\vec{\alpha}(s,\mu,y)_{3})^T.
\end{align*}
Then there are $\mu_1= \mu_1(t_0,k,h)>0$ and a $C^2$-function $\beta:(-s_0,s_0)\times (0,\mu_1) \to \rz^3$ 
depending on $t_0$, $k$, and $h$,
such that 
\begin{align*}
\beta(s,\mu) = -\mu^2 \frac12 \big(D^2k_\theta(0)\big)^{-1} \nabla \laplace k_\theta(0) + O(\mu^3), 
\end{align*}
as $\mu \to 0$ and
\begin{align*}
\hat{\alpha}(s,\mu,\beta(s,\mu)) =0 \text{ for all } (s,\mu) \in (-s_0,s_0)\times (0,\mu_1).
\end{align*}
Moreover, $\beta$ is unique in the sense that, 
if $y \in B_{\mu_1}(0)$ satisfies $\hat{\alpha}(s,\mu,y) =0$ 
for some $s \in (-s_0,s_0)$ and $0<\mu<\mu_1$, then $y= \beta(t,\mu)$. 
\end{lemma}
\begin{proof}
Lemma \ref{l:derivative_alpha_y} suggests to apply the implicit function theorem, 
but unfortunately $\vec{\alpha}$ may not be differentiable for $\mu=0$.
Instead we apply directly Banach's fixed-point 
theorem to the function 
\begin{align*}
F(s,\mu,y) := y + \big(D^2k_\theta(0)\big)^{-1}\hat{\alpha}(s,\mu,y)
\end{align*}
in $B_\delta(0)$, where $0<\delta<\dist(0,\supp(h_\theta))$ will be chosen later.\\
For $y \in B_\delta(0)$ we use the fact that $\nabla k_\theta(0)=0$ and get 
\begin{align}
\label{eq:17}
\bigg(\frac{1+t_0 k_\theta(0)}{1+t_0k_\theta(y)}\bigg)^{\frac54}= 1 +O(\delta^2).
\end{align}
Fix $y_1,y_2 \in B_\delta(0)$ and $(s,\mu )\in (-s_0,s_0)\times (0,\mu_0)$, then by Lemma \ref{l:derivative_alpha_y}
and (\ref{eq:17})
\begin{align*}
|F(s,&\mu,y_1)-F(s,\mu,y_2)|\\
&= \big|(y_1-y_2)+ \big(D^2k_\theta(0)\big)^{-1}  
\int_0^1 \frac{\rand \hat{\alpha}}{\rand y}(s,\mu,y_2+t(y_1-y_2)) (y_1-y_2)\di t \big|\\
&\le \big|(y_1-y_2)- \bigg(\int_0^1 \big(D^2k_\theta(0)\big)^{-1} D^2k_\theta(y_2+t(y_1-y_2)) \di t\bigg) (y_1-y_2)\big|\\
&\qquad +  O\Big(\delta^2+\supm_{y \in B_\delta(0)}|\nabla k_\theta(y)|+ \mu^\frac14\Big) |y_1-y_2|\\
&\le O\Big(\delta + \mu^\frac14\Big) |y_1-y_2|.  
\end{align*}
For $y \in B_\delta(0)$ we estimate using Lemma \ref{p:implicit}
\begin{align*}
|F(s,\mu,y)| &= \big|y + \big(D^2k_\theta(0)\big)^{-1}\big(\hat{\alpha}(s,\mu,y)\big)\big|\\
&\le \Big|y- \big(D^2k_\theta(0)\big)^{-1}\Big(\nabla k_\theta(y) + O(\delta^2+\mu^{2})\Big)\Big|\\
&\le \Big|y- \big(D^2k_\theta(0)\big)^{-1}\Big(D^2k_\theta(0)y+O(\delta^2+\mu^{2})\Big)\Big|\\
&\le O(\delta^2 +\mu^2).  
\end{align*}
Consequently, there is $\mu_1>0$ such that 
$F(s,\mu,\cdot)$ is a contraction in $B_{\mu_1}(0)$ for any $0 <\mu < \mu_1$ and $s\in [-s_0,s_0]$.
From Banach's fixed-point theorem we may define $\beta(s,\mu)$ to be the unique fixed-point of $F(s,\mu,\cdot)$
in $B_{\mu_1}(0)$. After shrinking $\mu_1$ if necessary we may apply 
Lemma \ref{l:derivative_alpha_y} and the usual implicit function theorem
to see that the function $\beta$ is twice differentiable for $\mu>0$.\\ 
To deduce the expansion for small $\mu$ we fix $\rho>0$ and 
\begin{align*}
y \in U_\rho:= B_\rho\Big(-\mu^2\frac12 \big(D^2k_\theta(0)\big)^{-1} \nabla \laplace k_\theta(0)\Big).
\end{align*}
Then, by Lemma \ref{p:implicit} and (\ref{eq:17})
\begin{align*}
\Big|F(s,&\mu,y)+\mu^2\frac12 \big(D^2k_\theta(0)\big)^{-1} \nabla \laplace k_\theta(0) \Big|\\
&\le \Big|y +\big(D^2k_\theta(0)\big)^{-1}
\Big(\hat{\alpha}(s,\mu,y)+\mu^2\frac12\nabla \laplace k_\theta(0)\Big)\Big|\\
&\le \Big|y + \big(D^2k_\theta(0)\big)^{-1}
\Big(-\nabla k_\theta(y)-\mu^2\frac12\big(\nabla \laplace k_\theta(y)-\nabla \laplace k_\theta(0)\big)\\
&\qquad \qquad \qquad \qquad \qquad   +O(\mu|y|^2+\mu^2|y|+\mu^{3})\Big)\Big|\\
&\le O\big(\rho^2 +\mu^2 \rho +\mu^{3}\big).
\end{align*}
Consequently, we may choose for small $0<\mu$ a radius $0<\rho= O(\mu^{3})$ such 
that $F$ maps $U_\rho \subset B_{\mu_1}(0)$ into itself.
Consequently, the unique fixed-point $\beta(s,\mu)$ must lie in this ball. This ends the proof.
\end{proof}  
Hence, to exclude or to construct blow-up sequences, which blow-up at a nondegenerate
critical point $\theta$ of $k$ with $\laplace k_\theta(0)=0$ it suffices to study 
$\a(s,\mu,\beta(s,\mu))_0$. 
\begin{lemma}
\label{l:expansion_tilde_alpha_0} 
Under the assumptions of Lemma \ref{l:implicit:curve} and $k\in C^5(S^3)$
we have
\begin{align*}
\big(\a(s,\mu,\beta(s,\mu))\big)_0
&= -t_0 \mu^3 (1+t_0 k(\theta))^{-\frac54} \frac{3^\frac34 4}{\pi \sqrt{5}} 
\Big(a_0(\theta) +s  \int_{\rz^3} h_\theta(x)|x|^{-6}\Big)\\
&\quad +t_0 \mu^4  
\frac{\pi 3^\frac34 \sqrt{5}}{30 (1+t_0 k(\theta))^{\frac94}} 
\Big(a_1(\theta)+t_0 a_2(\theta)\Big)
+ O(\mu^{4+\frac14}).
\end{align*}
\end{lemma}
\begin{proof}
In view of Lemma \ref{l:implicit:curve} and because $\nabla k_\theta(0)=0$
we may estimate functions of $\beta(s,\mu)$ and of $k(\beta(s,\mu))$ as follows
\begin{align}
\label{eq:105}
\begin{split}
&F(\beta)= F(0)-\mu^2 F'(0) \frac12\big(D^2k_\theta(0)\big)^{-1} \nabla \laplace k_\theta(0)
+ O(\mu^{3}),\\
&F(k(\beta))= F(k_\theta(0))+ O(\mu^{4}).     
\end{split}
\end{align}
To prove the claim of the lemma we expand $\a(s,\mu,\beta(s,\mu))_0$ according to Lemma \ref{p:implicit} 
and use \eqref{eq:105}. 
\end{proof}

\begin{lemma}
\label{l:d_gamma_d_eps}
Under the assumptions of Lemma \ref{l:implicit:curve} suppose 
\begin{align*}
a_0(\theta)=0 \text{ and } a_1(\theta)+t_0 a_2(\theta)>0,  
\end{align*}
and define 
\begin{align*}
\gamma(s,\mu):= -\frac{1}{t_0\mu^{3}} (1+t k(\theta))^{\frac54} \frac{\pi \sqrt{5}}{3^\frac34 4} 
\a(s,\mu,\beta(s,\mu))_0.
\end{align*}
Then as $\mu \to 0$
\begin{align}
\label{eq:18}
\frac{\rand \gamma(s,\mu)}{\rand s} &= \int h_\theta(x) |x|^{-6} + O(\mu),\\
\label{eq:19}
\frac{\rand \gamma(s,\mu)}{\rand \mu}&= -\frac{\pi^2}{24} (1+t_0 k(\theta))^{-1} (a_1(\theta)+t_0 a_2(\theta))
+O(\mu^{\frac14}).
\end{align}
\end{lemma}
\begin{proof}
As $\nabla k_\theta(0)=0$ we get from (\ref{eq:4}) that $\dist(0,\supp(h_\theta))>0$. As $\beta(s,\mu)=O(\mu^2)$
as $\mu \to 0$ we get that any term which depends only locally on $sh_\theta$ is independent of $s$ for small
$\mu>0$.\\ 
We have
\begin{align*}
\frac{d}{d s}\a(s,\mu,\beta(s,\mu))_0  
= \frac{\rand (\a)_0}{\rand s}\Big|_{(s,\mu,\beta(s,\mu))}
+ \frac{\rand (\a)_0}{\rand y}\Big|_{(s,\mu,\beta(s,\mu))} \frac{\rand \beta}{\rand s}\Big|_{(s,\mu)}. 
\end{align*}
The derivatives of $\a(\cdot)_0$ are given in Lemmas \ref{l:derivative_alpha_y}-\ref{l:derivative_alpha_s}.
To compute 
the derivative of $\beta$ we use the fact that 
$\hat{\a}(s,\mu,\beta(s,\mu)) \equiv 0$.
By (\ref{eq:105}) and Lemmas \ref{l:derivative_alpha_y}-\ref{l:derivative_alpha_s} we have
\begin{align*}
\frac{\rand \beta}{\rand s}\Big|_{(s,\mu)} &= 
-\bigg(\frac{\rand \hat{\a}}{\rand y}\Big|_{(s,\mu,\beta(s,\mu))}\bigg)^{-1}
\frac{\rand \hat{\a}}{\rand s}\Big|_{(s,\mu,\beta(s,\mu))}\\
&=\Big(\big(D^2k_\theta(0)\big)^{-1}+
O(\mu^{1+\frac14})\Big) \notag \\ 
&\quad \frac{3^\frac14 \sqrt{5}}{t_0 \mu \pi} \big(1+t_0 k(\theta)\big)^{\frac54}
\bigg[ \sum_{j=1}^3 \frac{\a_j(s,\mu,\beta)_i}{\rand s} + O(\mu^4)\bigg]_{i=1\dots 3}\\
&= O(\mu^{3}),
\end{align*}
where we used that $a_j(s,\mu,y)_i$ is independent of $s$ for small $|y|, \mu>0$. 
From Lemma \ref{l:derivative_alpha_y} we get
\begin{align*}
\frac{\rand (\a)_0}{\rand y}&\Big|_{(s,\mu,\beta(s,\mu))} \frac{\rand \beta}{\rand s}\Big|_{(s,\mu)}
= O(\mu^5).
\end{align*}
Furthermore, by Lemma \ref{l:derivative_alpha_s}
\begin{align*}
\frac{d \a(s,\mu,\beta)_0 }{ds}
&= \sum_{j=2}^4\frac{\rand \vec{\a_j}(s,\mu,\beta)_0}{\rand s} +O(\mu^{4+\frac14})\\
&= -t_0 \mu^3 \Big(1+t_0 k_\theta(0)\Big)^{-\frac{5}{4}} 
\frac{3^\frac34 4}{\pi \sqrt{5}}
\cpvint{}{} h_\theta(x) \frac{1}{|x|^{6}} +O(\mu^{4+\frac14}).
\end{align*}
The definition of $\gamma$, \eqref{eq:105}, and Lemma \ref{l:expansion_tilde_alpha_0}  
yield (\ref{eq:18}).\\
Concerning (\ref{eq:19}) we get
\begin{align*}
\frac{d}{d \mu}\a(s,\mu,\beta(s,\mu))_0  
= \frac{\rand (\a)_0}{\rand \mu}\Big|_{(s,\mu,\beta(s,\mu))}
+ \frac{\rand (\a)_0}{\rand y}\Big|_{(s,\mu,\beta(s,\mu))} \frac{\rand \beta}{\rand \mu}\Big|_{(s,\mu)}. 
\end{align*}
By (\ref{eq:105}) and Lemmas \ref{l:derivative_alpha_y}-\ref{l:derivative_alpha_mu} we have
\begin{align*}
\frac{\rand \beta}{\rand \mu}\Big|_{(s,\mu)} &= 
-\bigg(\frac{\rand \hat{\a}}{\rand y}\Big|_{(s,\mu,\beta(s,\mu))}\bigg)^{-1}
\frac{\rand \hat{\a}}{\rand \mu}\Big|_{(s,\mu,\beta(s,\mu))}\\
&=\Big(\big(D^2k_\theta(0)\big)^{-1}+
O(\mu^{1+\frac14})\Big) \notag \\ 
&\quad \frac{3^\frac14 \sqrt{5}}{t_0 \mu \pi} \big(1+t_0 k(\theta)\big)^{\frac54}
\bigg[ \sum_{j=1}^3 \frac{\a_j(s,\mu,\beta)_i}{\rand \mu} + O(\mu^3)\bigg]_{i=1\dots 3}\\
&= \Big(\big(D^2k_\theta(0)\big)^{-1}+
O(\mu^{1+\frac14})\Big) \notag \\ 
&\quad \bigg(
-\mu\nabla\laplace k_\theta(0)
+ \frac{3^\frac14 \sqrt{5}}{t_0 \mu \pi}\bigg[ \frac{1}{\mu}\sum_{j=1}^3 \a_j(s,\mu,\beta)_i +  O(\mu^3)\bigg]_{i=1\dots 3}
\bigg)\\
&=-\mu \big(D^2k_\theta(0)\big)^{-1}\nabla\laplace k_\theta(0)+ O(\mu^2). 
\end{align*}
Hence, by Lemma \ref{l:derivative_alpha_y} and \ref{l:derivative_alpha_mu}
\begin{align*}
&\frac{\rand (\a)_0}{\rand y}\Big|_{(s,\mu,\beta(s,\mu))} \frac{\rand \beta}{\rand \mu}\Big|_{(s,\mu)}\\
&\qquad =\frac{t_0 \mu^3 \pi}{3^\frac14 \sqrt{5}} \Big(1+t_0k_\theta(0)\Big)^{-\frac54}
\nabla \laplace k_\theta(0) \big(D^2k_\theta(0)\big)^{-1}\nabla\laplace k_\theta(0)   
+ O\big(\mu^{3+\frac14}\big), 
\end{align*}
and
\begin{align*}
 \frac{\rand \gamma(s,\mu)}{\rand \mu}&=  
\frac{3}{t_0\mu^{4}} (1+t k(\theta))^{\frac54} \frac{\pi \sqrt{5}}{3^\frac34 4} 
\a(s,\mu,\beta(s,\mu))_0\\ 
&\quad -\frac{1}{t_0\mu^{3}} (1+t k(\theta))^{\frac54} \frac{\pi \sqrt{5}}{3^\frac34 4}
\sum_{j=2}^4 \frac{\rand \alpha_j(s,\mu,\beta)_0}{\rand \mu}\\ 
&\quad -\frac{\pi^2}{12} \nabla \laplace k_\theta(0) \big(D^2k_\theta(0)\big)^{-1}\nabla\laplace k_\theta(0)
+O\big(\mu^{3+\frac14}\big)\\
&= \frac{3}{t_0\mu^{4}} (1+t k(\theta))^{\frac54} \frac{\pi \sqrt{5}}{3^\frac34 4} 
\sum_{j=2}^4 \alpha_j(s,\mu,\beta)_0
-\frac{(1+t k(\theta))^{\frac54}}{t_0\mu^{4}} \frac{\pi \sqrt{5}}{3^\frac34 4}\\
&\quad \cdot \Big({3}\sum_{j=2}^4 \alpha_j(s,\mu,\beta)_0 -\alpha_2(s,\mu,\beta)_0+\alpha_4(s,\mu,\beta)_0\Big)\\
&\quad -\frac{\pi^2}{12} \nabla \laplace k_\theta(0) \big(D^2k_\theta(0)\big)^{-1}\nabla\laplace k_\theta(0)
+O\big(\mu^{\frac14}\big)\\
&= -\frac{1}{t_0\mu^{4}} (1+t k(\theta))^{\frac54} \frac{\pi \sqrt{5}}{3^\frac34 4}
\Big(-\alpha_2(s,\mu,\beta)_0+\alpha_4(s,\mu,\beta)_0\Big)\\
&\quad -\frac{\pi^2}{12} \nabla \laplace k_\theta(0) \big(D^2k_\theta(0)\big)^{-1}\nabla\laplace k_\theta(0)
+O\big(\mu^{\frac14}\big) 
\end{align*}
If we use \eqref{eq:105} and the expansion in Lemma \ref{p:implicit} we find
\begin{align*}
&-\frac{1}{t_0\mu^{4}} (1+t k(\theta))^{\frac54} \frac{\pi \sqrt{5}}{3^\frac34 4}
\Big(-\alpha_2(s,\mu,\beta)_0+\alpha_4(s,\mu,\beta)_0\Big)\\
&\quad = \frac{\pi^2}{24} \nabla \laplace k_\theta(0) \big(D^2k_\theta(0)\big)^{-1}\nabla\laplace k_\theta(0)\\
&\quad \quad -\frac{\pi^2}{24} \laplace^2 k_\theta(0)+\frac{t_0}{1+t_0 k(\theta)} \frac{5\pi}{64}
\intg_{\rand B_1(0)} \big|D^2 k_\theta(0)(x)^2\big|^2 \di Sx.
\end{align*}
Summing up yields the claim of the lemma.
\end{proof}

\begin{lemma}
\label{l:blow_up_curves}
Under the assumptions of Lemma \ref{p:implicit}
we define $M_* \subset S^3$ by
\begin{align}
M_*:= \{\theta \in \text{Crit}(k) \where \laplace k_\theta(0)=0=a_0(\theta), a_1(\theta)+t_0 a_2(\theta)>0\}.
\end{align}
Then there is $\delta>0$ such that for any $\theta \in M_*$
there exists a unique $C^1$-curve
\begin{align*}
\{0<\mu<\delta\} \ni \mu \mapsto (s^\theta(\mu),\phi^\theta(\mu,\cdot)) 
\in (0,\delta)\times C^{2,\a}(S^3), 
\end{align*}
such that as $\mu \to 0$
\begin{align*}
s^\theta(\mu) = \mu \frac{\pi^2}{24}\Big(\int h_\theta(x) |x|^{-6}\Big)^{-1} 
\frac{a_1(\theta)+t_0 a_2(\theta)}{1+t_0 k(\theta)}
+O(\mu^{1+\frac14}),\\
\frac{\rand s^\theta}{\rand \mu}(\mu) = \frac{\pi^2}{24}\Big(\int h_\theta(x) |x|^{-6}\Big)^{-1} 
\frac{a_1(\theta)+t_0 a_2(\theta)}{1+t_0 k(\theta)}
+O(\mu^{\frac14}).  
\end{align*}
and $\phi^\theta(\mu,\cdot)$ solves \eqref{eq:eq1} for $s=s^\theta(\mu)$ and blows up like
\begin{align*}
\|{\mathcal R}_\theta(\phi^\theta(\mu,x))-(1+t_0 k(\theta))^{-\frac14} 
z_{\mu,0}(x)\|_{\Did\cap C^2(B_1(0))} = O(\mu^2).  
\end{align*} 
The curves are unique, in the sense that, 
if $(s_i,\phi_i) \in (0,\delta) \times C^{2,\a}(S^3)$ blow
up at some $\theta \in S^3$ then $\theta \in M_*$ and there is a sequence of positive numbers 
$(\mu_i)$ converging to zero such that 
$(s_i,\phi_i) = (s^\theta(\mu_i),\phi^\theta(\mu_i,\cdot))$
for all but finitely many $i \in \nz$.
\end{lemma}
\begin{proof}
We fix $\theta \in M_*$.
To construct $s^\theta(\mu)$ we proceed as in Lemma \ref{l:implicit:curve} and 
use Banach's fixed-point theorem applied to
\begin{align*}
F_2(s,\mu) := s- \Big(\int h_\theta(x) |x|^{-6}\Big)^{-1} \gamma(s,\mu).  
\end{align*}
Since we know the expansion of $\gamma$ and $\frac{\rand \gamma}{\rand s}$ as $\mu \to 0$
it is easy to see that $F_2(\cdot,\mu)$ is a contraction in
\begin{align*}
B_r\Bigg(\mu \frac{\pi^2}{24}\Big(\int h_\theta(x) |x|^{-6}\Big)^{-1} 
\frac{a_1(\theta)+t_0 a_2(\theta)}{1+t_0 k(\theta)}\Bigg)  
\end{align*}
for any $0<\rm{const}\, \mu^{1+\frac14}\le r\le r_1$ and the existence part of the claim follows from that.
The differentiability of $s$ with respect to $\mu$ follows from Lemma \ref{l:d_gamma_d_eps} and the usual
implicit function theorem.\\
Assume $(s_i,\phi_i)$ blow up at some $\theta \in S^3$. Then we apply Lemma \ref{l:blow_up} and find in
stereographic coordinates $\stereo_\theta$ sequences $y_i \to 0$, $\mu_i \to 0$ such that
\begin{align*}
{\mathcal R}_\theta(\phi_i)(x)-\big(1+t_0 (k(\theta)+s_ih(\theta)\big)^{-\frac14} 
z_{\mu_i,y_i}(x)   
\end{align*}
is orthogonal to $T_{\mu_i,y_i}Z$ and converges to $0$ as $i \to \infty$. Consequently,
if we set 
\begin{align*}
w(i):= {\mathcal R}_\theta(\phi_i)-z_{\mu_i,y_i}  
\end{align*}
we find as $z_{\mu_i,y_i}$ is orthogonal to $T_{\mu_i,y_i}Z$, 
\begin{align*}
w(i) \text{ is orthogonal to } T_{\mu_i,y_i}Z \ \text{ and } \ w(i)-w_0(s_i,\mu_i,y_i) \to_{i \to \infty}  0, 
\end{align*}
where $w_0$ is defined in Lemma \ref{p:implicit}. Moreover
\begin{align*}
0=f_{t_0,s}'({\mathcal R}_\theta(\phi_i)) = f_{t_0,s}'(z_{\mu,y}+w(i)). 
\end{align*}
The uniqueness part of Lemma \ref{p:implicit} shows for large $i$
\begin{align*}
w(i)= w(s_i,\mu_i,y_i) \text{ and } \vec{\a}(s_i,\mu_i,y_i)=0.  
\end{align*}
As $\mu_i \to 0$ the expansion of $\vec{\a}$ of order $\mu$ and $\mu^{2}$ in Lemma \ref{p:implicit} shows
\begin{align*}
\nabla k_\theta(0) =0 \text{ and } \laplace k_\theta(0)=0.  
\end{align*}
From Lemma \ref{l:implicit:curve} we infer that
\begin{align*}
y_i = \beta(s_i,\mu_i)  
\end{align*}
and the expansion in Lemma \ref{l:expansion_tilde_alpha_0} gives
\begin{align*}
0 &= -t_0 \mu_i^3 (1+t_0 k(\theta))^{-\frac54} \frac{3^\frac34 4}{\pi \sqrt{5}} 
\Big(a_0(\theta) +s_i  \int_{\rz^3} h_\theta(x)|x|^{-6}\Big)\\
&\quad +t_0 \mu_i^4 (1+t_0 k(\theta))^{-\frac94} \frac{\pi 3^\frac34 \sqrt{5}}{30} \Big(a_1(\theta)+t_0 a_2(\theta)\Big)
+ O(\mu_i^{4+\frac14}).  
\end{align*}
Consequently 
\begin{align*}
\Big(a_0(\theta) +s_i  \int_{\rz^3} h_\theta(x)|x|^{-6}\Big) \to 0 \text{ as } i \to \infty, 
\end{align*}
and from the choice of $h$, assuming $0<\delta<s_0$, 
we deduce that $a_0(\theta)=0$. Hence
\begin{align*}
s_i \int_{\rz^3} h_\theta(x)|x|^{-6}
=
\mu_i (1+t_0 k(\theta))^{-1}
\frac{\pi^2}{24} \Big(a_1(\theta)+t_0 a_2(\theta)\Big)
+ O(\mu_i^{1+\frac14}).
\end{align*}
Thus, $a_1(\theta)+t_0 a_2(\theta)$ has to be positive, which shows $\theta \in M_*$, and for large $i$
\begin{align*}
s_i \in B_{r_1} \Bigg(\mu_i \frac{\pi^2}{24}\Big(\int h_\theta(x) |x|^{-6}\Big)^{-1} 
\frac{a_1(\theta)+t_0 a_2(\theta)}{1+t_0 k(\theta)}\Bigg). 
\end{align*}
The uniqueness of the fixed point implies $s_i=s_i^\theta(\mu_i)$ and the claim follows.
\end{proof}


\section{The Leray-Schauder degree}
\label{s:degree}
From Section \ref{sec:introduction} we know that
the degree 
$\deg(Id-L_s,\mathcal{B}_{C_{\delta}},0)$ of the problem (\ref{eq:10}) is independent of $s \in [\delta,s_0]$ and equals
\begin{align*}
\deg(Id-L_s,\mathcal{B}_{C_{\delta}},0)= -\Big(1+ \sum_{\theta \in Crit_-(k+sh)} (-1)^{\ind(k,\theta)}\Big), 
\end{align*}
where the set $Crit_-(k+sh)$ is independent of $s$ and given by
\begin{align*}
Crit_-(k+sh)
&= \Big\{ \theta \in \text{Crit}(k) \where \laplace k(\theta)<0 \text{ or }\\ 
&\qquad \big(\laplace k(\theta)=0 \text{ and } a_0(\theta)<0\big)\Big\}.
\end{align*}
By Lemma \ref{l:blow_up_curves} and the a priori estimate for $s=0$ the set of functions
\begin{align*}
L_b := \big\{ \phi \text{ solves } &(\ref{eq:10}) \text{ for some } s \in [0,s_0],\,\\
&\phi \not\in \cup_{\theta \in M_*} 
\{\phi^\theta(s^\theta(\mu),\cdot)\where 0<\mu<\delta\}\big\}   
\end{align*}
is uniformly bounded from above and by standard elliptic regularity also in $C^{2,\alpha}(S^3)$. 
By Sobolev's and Harnack's inequality this gives a uniform lower bound,
thus there is $C_1>0$ such that
$
L_b \subset  \mathcal{B}_{C_1}. 
$\\
Again from Lemma \ref{l:blow_up_curves} and since $\frac{\rand s^\theta}{\rand \mu}$ is uniformly positive,
there is $s_1>0$ small, such that for any $0<s\le s_1$ and any $\theta \in M_*$ 
there exists exactly one $\mu^\theta(s) \in (0,\delta)$
satisfying
\begin{align*}
s^\theta\big(\mu^\theta(s)\big)=s.  
\end{align*}
Moreover, we may assume, shrinking $s_1$
\begin{align*}
\|\phi^\theta(\mu^\theta(s),\cdot)\|_\infty \ge 2C_1 \; \forall \theta \in M_*,\\
\|\phi^{\theta_1}(\mu^{\theta_1}(s),\cdot)-\phi^{\theta_2}(\mu^{\theta_1},\cdot)\|_\infty \ge C_1
\; \forall \theta_1\neq \theta_2 \text{ in } M_*.  
\end{align*}
Hence, there are two types of solutions to (\ref{eq:10}) as $s \to 0^+$: the solutions in
$L_b\subset \mathcal{B}_{C_1}$ remain uniformly bounded as $s\to 0^+$ 
and the solutions $\{\phi^\theta(\mu^\theta(s),\cdot)\where \theta \in M_*\}$ that blow up as $s\to 0^+$ and are 
uniformly isolated for each fixed small $s>0$.
Consequently, using the additivity of the degree, we find for any $0<s\le s_1$
\begin{align*}
\deg(Id-L_s,&\mathcal{B}_{C_{s}},0)\\
&= \deg(Id-L_s,\mathcal{B}_{C_{1}},0) +\sum_{\theta \in M_*} \deg_{loc}(Id-L_s,\phi^\theta(\mu^\theta(s),\cdot))\\
&= \deg(Id-L_0,\mathcal{B}_{C_{1}},0) +\sum_{\theta \in M_*} \deg_{loc}(Id-L_s,\phi^\theta(\mu^\theta(s),\cdot)).
\end{align*}
Together with (\ref{eq:16}) we get for any $0<s\le s_1$
\begin{align*}
\deg&(Id-L_0,\mathcal{B}_{C_{1}},0) \\
&= 
-\Big(1+ \sum_{\theta \in Crit_-(k+sh)} (-1)^{\ind(k,\theta)}\Big)
-\sum_{\theta \in M_*} \deg_{loc}(Id-L_s,\phi^\theta(\mu^\theta(s),\cdot)).
\end{align*}
It remains to compute the local degree $\deg_{loc}(Id-L_s,\phi^\theta(\mu^\theta(s),\cdot))$ for any
$\theta \in M_*$. 
We use the transformation ${\mathcal R}_\theta$ in (\ref{eq:116}) to define
the weighted space 
\begin{align*}
  C^{2}(\rz^{3},{\mathcal R}_\theta)
&:= \big\{ u \in C^{2}(\rz^{3}) \where u \in {\mathcal R}_\theta\big(C^{2}(S^{3})\big)\big\},\\
\|u\|_{C^{2}(\rz^{3},{\mathcal R}_\theta)} &:= \|({\mathcal R}_\theta)^{-1}(u)\|_{C^{2}(S^{3})}.
\end{align*}
Note that $C^{2}(\rz^{3},{\mathcal R}_\theta)\embed \Did$, because ${\mathcal R}_\theta$ is an 
isomorphism between $H^{1,2}(S^{3})$ and $\Did$. Using ${\mathcal R}_\theta$ we obtain
\begin{align*}
\deg_{loc}(Id-L_s,\phi^\theta(\mu^\theta(s),\cdot))
&= \deg_{loc}(Id- {\mathcal R}_\theta L_s ({\mathcal R}_\theta)^{-1} ,u_{\theta,s})\\
&= \deg_{loc}\Big(f_{t_0,s}',u_{\theta,s}\Big),
\end{align*}
where $u_{\theta,s}= {\mathcal R}_\theta \big(\phi^\theta(\mu^\theta(s),\cdot)\big)  \in C^{2}(\rz^{3},{\mathcal R}_\theta)$. 
Note that by duality we consider $f_{t_0,s}'$ as a map from the Hilbert space $D^{1,2}(\rz^3)$ into itself.
\begin{lemma}
\label{l:deg_blow_up_solutions}
Under the assumptions of Lemma \ref{l:blow_up_curves} there holds for $0<s\le s_1$ 
\begin{align*}
\sum_{\substack{\theta \in M_*}}  
\deg_{loc}\Big(f_{t_0,s}',u_{\theta,s}\Big)
= \sum_{\substack{\theta \in M_*}} (-1)^{\ind(k,\theta)}.
\end{align*}
\end{lemma}
\begin{proof}
Fix $\theta \in M_*$. The solution $u_{\theta,s} \in C^{2}(\rz^{3},{\mathcal R}_\theta)$
is given in notation of Lemmas \ref{p:implicit} and \ref{l:implicit:curve} by  
\begin{align*}
u_{\theta,s} = z_{\mu^\theta(s),y^\theta(s)}+w\big(s,\mu^\theta(s),y^\theta(s)\big),
\end{align*}
where $y^{\theta}(s) = \beta^{\theta}(s,\mu^{\theta}(s))$. 
$(\mu^\theta(s),y^\theta(s))$ is the only zero of $\vec{\alpha}(s,\cdot,\cdot)$ 
for $\mu$ and $|y|$ bounded above by a small fixed constant.
As $y^{\theta}(s)=O(s^{2})$ we may replace $y^{\theta}$ by $0$ 
(in various expressions below) and get an addition $O(s^{2})$-error.\\
We drop the $s$-dependence of $\mu^\theta$ and $y^\theta$ in the notation 
when there is no possibility of confusion.
Moreover by Lemma \ref{l:blow_up_curves} we have $s^\theta(\mu)\thicksim \mu$ and we 
may estimate the errors in terms of $s$.\\ 
As seen above by Lemma \ref{l:blow_up_curves} 
the solution $u_{\theta,s}$ remains uniform isolated in $C^{2}(\rz^{3},{\mathcal R}_\theta)$ as well as
in $\Did$
for $s \in (0,s_1]$. From (\ref{eq:13}) and regularity results \cite{Luckhaus79,BrezisKato79}
we infer that $w(s,\mu,y) \in C^{2}(\rz^{3},{\mathcal R}_\theta)$
depends continuously on $(s,\mu,y)$.\\
To compute the local degree, we first show that $f_{t_0,s}''(u_{\theta,s})$ is nondegenerate.
To this end we let 
\begin{align*}
\varphi(s,\theta)_0 &:= 
\mu^{\theta} c_\xi^{-1} 
\frac{\rand}{\rand\mu} \big(z_{\mu,\beta(s,\mu)}+w(s,\mu,\beta(s,\mu))\big)\eval_{\mu^{\theta}},\\
\varphi(s,\theta)_i &:=
\mu^{\theta} c_\xi^{-1} 
\frac{\rand}{\rand y_i}\big(z_{\mu^{\theta},y}+w(s,\mu^{\theta},y)\big)\eval_{y^{\theta}}, \quad i=1\dots 3.
\end{align*}
The derivatives of $\beta$ and $w$ with respect to $\mu$ are computed in \cite[App. A]{schneider04}
the derivatives of $w$ with respect to $y_i$ are given in  \cite[Lem. 5.1]{schneider04_2}. We have
\begin{align*}
|\frac{\rand\beta}{\rand \mu}(s,\mu^{\theta})| &= O(s),\\
\|\frac{\rand w}{\rand \mu}(s,\mu^{\theta},y^{\theta})-\frac{\rand w_0}{\rand \mu}(s,\mu^{\theta},y^{\theta})\|_{\Did}&= O(s)\\
\|\frac{\rand w}{\rand y_i}(s,\mu^{\theta},y^{\theta})-\frac{\rand w_0}{\rand y_i}(s,\mu^{\theta},y^{\theta})\|_{\Did}&= O(s).  
\end{align*}
Therefore we get
\begin{align*}
u_{\theta,s} &= (1+t_0 k_\theta(0))^{-\frac14}z_{\mu^{\theta},y^{\theta}} +O(s^{2})_{\Did},\\
\varphi(s,\theta)_0
&= (1+t_0 k_\theta(0))^{-\frac14}(\dot\xi_{\mu^{\theta},y^{\theta}})_0 +O(s^{2})_{\Did},\\
\varphi(s,\theta)_i
&= (1+t_0 k_\theta(0))^{-\frac14}(\dot\xi_{\mu^{\theta},y^{\theta}})_i +O(s^{2})_{\Did},
\end{align*}
By Lemma \ref{l:blow_up_curves} and Lemma \ref{p:implicit} we find
\begin{align*}
f'_{t_0,s}(z_{\mu,\beta(s,\mu)}+w(s,\mu,\beta(s,\mu))) = \alpha(s,\mu,\beta(s,\mu))_0 (\dot\xi_{\mu,\beta(s,\mu)})_0.
\end{align*}
Differentiating with respect to $\mu$ by Lemma \ref{l:d_gamma_d_eps} leads to
\begin{align*}
f''_{t_0,s}&(u_{\theta,s}) 
(\mu^{\theta})^{-1} c_\xi \varphi(s,\theta)_0\\
&= \Big(t_0(\mu^{\theta})^{3} (1+t k_\theta(0))^{-\frac94} \frac{3^\frac34 }{6\sqrt{5}} 
   (a_1(\theta)+t_0 a_2(\theta)) + O(s^{3+\frac14})\Big)
(\dot\xi_{\mu^{\theta},y^{\theta}})_0.
\end{align*}
Moreover, differentiating
\begin{align*}
f'_{t_0,s}(z_{\mu,y}+w(s,\mu,y)) = \sum_{i=0}^{3}\alpha(s,\mu,y)_i (\dot\xi_{\mu,y})_i  
\end{align*}
with respect to $y_j$ we get from Lemma \ref{l:derivative_alpha_y}
\begin{align*}
f''_{t_0,s}(u_{\theta,s})\frac{c_\xi}{\mu^{\theta}} \varphi(s,\theta)_j
&= 
-\frac{t_0 \mu^{\theta} \pi  
\sum_{i=1}^{3}
\Big(\frac{\rand^2 k_\theta }{\rand x_i \rand x_j}(0)+O(s)\Big) 
(\dot\xi_{\mu^{\theta},y^{\theta}})_i}
{3^\frac14 \sqrt{5} (1+t_0 k_\theta(0))^{\frac54}}
\\
&\quad +O(s^{2}) (\dot\xi_{\mu^{\theta},y^{\theta}})_0.
\end{align*}
Orthogonal to $T_{z_{\mu^{\theta},y^{\theta}}}Z$ we use
\begin{align}
\notag
f_{t_0,s}''(u_{s,\theta}) &= f_0''(z_{\mu^{\theta},y^{\theta}})
+ O(\|w(s,\mu^{\theta},y^{\theta})-w_0(s,\mu^{\theta},y^{\theta})\|)_{{\mathcal L}(\Did)}\\
\notag
&\quad -\frac{5t_0}{(1+t_0k_\theta(y^{\theta}))} 
\int_{\rz^{3}}((k_\theta+sh_\theta)(x)-k_\theta(y)) (z_{\mu^{\theta},y^{\theta}})^{4} \cdot \cdot \, dx\\
\label{eq:1}
&= f_0''(z_{\mu^{\theta},y^{\theta}}) + O(\mu^{\theta})_{{\mathcal L}(\Did)}.   
\end{align}
The $O(\mu)$-estimates are given in \cite[Lem. 4.1]{schneider04_2} or can be obtained using 
H\"older's and Sobolev's inequality and the fact that $k_\theta(x)-k_\theta(y)$ is bounded
in $\rz^{3}$ and of order $O(|x-y|)$ for $|x-y|<<1$.\\
To obtain a contradiction assume there is a function $v \in C^2(\rz^3,{\mathcal R}_\theta) \setminus \{0\}$ with
$f_{t_0,s}''(u_{\theta,s})v=0$. We may assume $\|v\|_{\Did}=1$. Then by (\ref{eq:1})
\begin{align*}
O(s) = \|f_0''(z_{\mu^{\theta},y^{\theta}})v\|_{\Did} 
\ge c \|Proj_{T_{z_{\mu^{\theta},y^{\theta}}}Z^\perp}v\|_{\Did},
\end{align*}
because $f_0''(z_{\mu^{\theta},y^{\theta}})$ is an isomorphism of 
$T_{z_{\mu^{\theta},y^{\theta}}}Z^\perp$. Moreover,
\begin{align*}
0 &= f''_{t_0,s}(u_{\theta,s}) c_\xi \varphi(s,\theta)_0 v\\
&= \Big( \frac{3^\frac34 t_0(\mu^{\theta})^{4} }{6\sqrt{5} (1+t k_\theta(0))^{\frac94}} 
   (a_1(\theta)+t_0 a_2(\theta)) + O(s^{3+\frac14})\Big)
\langle (\dot\xi_{\mu^{\theta},y^{\theta}})_0,v\rangle_{\Did}, 
\end{align*}
and 
\begin{align*}
\vec{0} &= \Big(f''_{t_0,s}(u_{\theta,s}) c_\xi \varphi(s,\theta)_j v\Big)_j\\
&= -\frac{t_0 (\mu^{\theta})^2 \pi}{3^\frac14 \sqrt{5} (1+t_0 k_\theta(0))^{-\frac54}}
\Big(D^2k_\theta(0)+O(s)\Big) \big(\langle(\dot\xi_{\mu^{\theta},y^{\theta}})_i,v\rangle_{\Did}\big)_i \\
&\quad +O(s^{3}) \langle(\dot\xi_{\mu^{\theta},y^{\theta}})_0,v\rangle_{\Did}.
\end{align*}
Since $D^2k_\theta(0)$ is invertible, we see that $Proj_{T_{z_{\mu^{\theta},y^{\theta}}}Z}v=0$,
contradicting the fact that $\|v\|_{\Did}=1$. Since $f''_{t_0,s}(u_{\theta,s})$ is of the form $id-compact$
in $C^2(\rz^3,{\mathcal R}_\theta)$ (as well as in $\Did$) we get
\begin{align*}
\|f''_{t_0,s}(u_{\theta,s})v\|_{C^2(\rz^3,{\mathcal R}_\theta)} \ge c \|v\|_{C^2(\rz^3,{\mathcal R}_\theta)}.   
\end{align*}
For $f_{t_0,s}'(u)=f_{t_0,s}''(u_{\theta,s})(u-u_{\theta,s})+O(\|u-u_{\theta,s}\|^2_{C^2(\rz^3,{\mathcal R}_\theta)})$,
\begin{align*}
\deg_{loc}\Big(f_{t_0,s}',u_{\theta,s}\Big) = \deg_{loc}\Big(f_{t_0,s}''(u_{\theta,s}),0\Big).  
\end{align*}
To compute $\deg_{loc}\Big(f_{t_0,s}''(u_{\theta,s}),0\Big)$ we consider the finite dimensional 
spaces (see (\ref{eq:5}))
\begin{align*}
X_{n,s} := &\langle u_{\theta,s}\rangle 
\oplus 
\langle 
\varphi(s,\theta)_0
\rangle
\oplus
\langle 
\varphi(s,\theta)_i:\, 1 \le i\le 3
\rangle\\
&\oplus
\langle \Phi^{\mu^{\theta},y^{\theta}}_{i,j,l}:\, i,j\in \nz_0,\, 2\le i+j\le n ,\, 1\le l\le  c_i \rangle.    
\end{align*}
The functions, spanning $X_{n,s}$, are a basis, as they are
orthogonal in $\Did$ up to an $O(s^{2})$-error.  The linear operator $Proj_{X_{n,s}}f_{t_0,s}''(u_{\theta,s})$ 
restricted to $X_{n,s}$ is given by, up to a multiplication of the elements in the diagonal by
positive constants
\begin{align*}
\begin{pmatrix}
-4& 0       & 0& 0\\
0& \mu^{4}(a_1(\theta)+t_0a_2(\theta))& 0&0\\
0& 0 & -\mu^{2}D^{2}k_\theta(0) & 0\\
0& 0 & 0 & f_0''(z_{\mu^{\theta},y^{\theta}})\eval_{\langle \Phi^{\mu^{\theta},y^{\theta}}_{i,j,l}\rangle}
\end{pmatrix}\\
+
\begin{pmatrix}
O(\mu) & O(\mu)       & O(\mu)& O(\mu)\\
O(\mu^{6})& O(\mu^{4+\frac14}) & O(\mu^{6}) &O(\mu^{6})\\
O(\mu^{4}) & O(\mu^{3}) 
& O(\mu^{3})  & O(\mu^{4})\\
O(\mu) & O(\mu) & O(\mu) & O(\mu)
\end{pmatrix}.
\end{align*}
Thus, we find for large $n$ and small $s$
\begin{align*}
\deg_{loc}\Big(f_{t_0,s}''(u_{\theta,s}),0\Big) &=
\sgn \det\big(Proj_{X_{n,s}}f_{t_0,s}''(u_{\theta,s})\big)\\
&= \sgn \det(D^{2}k_\theta(0)) = (-1)^{\ind(k,\theta)},
\end{align*}
which proofs the claim. 
\end{proof}
\begin{remark}
From the proof of Lemma \ref{l:deg_blow_up_solutions} we see
that $f_{t_0,s}''(u_{\theta,s})$ is nondegenerate and the Morse-Index
of $u_{\theta,s}$, i.e. the number of negative eigenvalues of $f_{t_0,s}''(u_{\theta,s})$,
is given by
\begin{align*}
\ind(f_{t_0,s},u_{\theta,s}) = 1+\ind(-k,\theta) = 4-\ind(k,\theta).  
\end{align*}
\end{remark}



\end{document}